\documentclass[]{article} 
\pdfoutput=1
\usepackage{authblk}
\usepackage{graphicx}
\usepackage{natbib}
\setcitestyle{aysep={}} 
\usepackage{amsmath}
\usepackage{amsthm}
\usepackage{amssymb}
\usepackage{enumerate}
\usepackage{enumitem}
\usepackage{subfig}
\usepackage{bbm}
\usepackage{todonotes}
\usepackage[hidelinks]{hyperref}

\newtheoremstyle{mystyle}
{10pt} 
{10pt} 
{\normalfont} 
{} 
{\bfseries} 
{.} 
{.5em} 
{} 

\theoremstyle{mystyle}

\newtheorem{theorem}{Theorem}

\newtheorem{proposition}{Proposition}

\theoremstyle{break}

\newcommand{\transp}[1]{#1^\mathsf{T}}
\newcommand{\p}[1]{\left(#1\right)}

\newcommand{\pow}[0]{\kern 0.12em}
\newcommand\ind{\mathbbm{1}}%

\providecommand{\keywords}[1]{\textbf{Keywords:} #1}

\setenumerate{leftmargin=1.7 cm}

\begin{document}

\title{SIR epidemics and vaccination on random graphs with clustering}

\author[1]{Carolina Fransson}
\author[2]{Pieter Trapman}

\affil[1]{\footnotesize Department of Mathematics, Stockholm University, 106 91 Stockholm, Sweden\\  \href{carolina.fransson@math.su.se}{carolina.fransson@math.su.se}} 
\affil[2]{\footnotesize Department of Mathematics, Stockholm University, 106 91 Stockholm, Sweden\\  \href{ptrapman@math.su.se}{ptrapman@math.su.se}}


%







\maketitle


\begin{abstract}

In this paper we consider SIR (Susceptible $\to$ Infectious $\to$ Recovered) epidemics 
on random graphs with clustering. 
 To incorporate group structure of the underlying social network, we use a generalized version of the configuration model in which each node is a member of a specified number of triangles. 
 SIR epidemics on this type of graph have earlier been investigated under the assumption of homogeneous infectivity and also under the assumption of  Poisson transmission and recovery rates.
 
We extend 
known results from literature by relaxing the assumption of homogeneous infectivity. 
An important special case of the epidemic model analyzed in this paper is epidemics in continuous time with  arbitrary infectious period distribution.
We use branching process approximations of the spread of the disease to provide 
expressions for the 
 basic reproduction number $R_0$, the  probability of a major outbreak and the expected final size. 
In addition, the impact of random vaccination with a perfect vaccine on the final outcome of the epidemic is investigated. 
We find that, for this particular model, $R_0$ equals the perfect vaccine-associated reproduction number. 

Generalizations to groups larger than three are discussed briefly. 

\end{abstract}


\section{Introduction}

One of the most important factors that determine the fate of an outbreak of an infectious disease is the contact pattern of individuals in the population. The frequency and duration of the contacts between individuals typically depend on the nature of their relationship. 
For this reason, recent interest has focused on the impact of the underlying social network on the spread of the disease. 
The social network is typically represented by a random graph  \citep{Newman2002}, in which the nodes or vertices represent individuals and the edges represent social contacts between the individuals. Two nodes that share an edge are called ``neighbors''.

A popular choice when generating 
random graphs with a specified degree distribution is the configuration model (CM). It was introduced by 
\citet{config} for the special case where the degree distribution is degenerate (i.e. every node of the graph has the same degree)
and extended to more general degree distributions by
\citet{MolloyReed95, MolloyReed98}.  
There is a vast literature on epidemics on configuration model graphs (see e.g.\ \citet{mathscientist, CM_vacc, CMJanson, barbour2013, front}).

An important feature of the configuration model is that, under mild regularity conditions on the degrees, this type of graph is asymptotically unclustered. That is to say, 
it contains virtually no groups and short circuits. 
Real world networks do, however, typically exhibit clustering \citep{Newman2003},  
and there are a number of 
 graph models that do allow for group structure 
 \citep{sparseclust, intersection, clustnewman2002}. 
Epidemics on graphs with group structure were studied by 
 \citet{TRAPMAN2007, backwardsimple,   backwardhousehold, backwardadv, Contagions, tunableclustering}.

In this paper, we use a generalized version of the configuration model to incorporate clustering of the social network
 in the analysis of the spread of an infectious disease. 
The configuration model with clustering (CMC)  was independently introduced by \citet{model_miller} and \citet{newmanconf}.
It is an extension of the CM in the sense that, for each node $u$, in addition to the degree of $u$ one also specifies the number of pairs of neighbors of $u$ that are in turn neighbors of each others.
In other words, one specifies the number of triangles (with non-overlapping edges) of which $u$ is a member (see section \ref{sec:model} for a precise definition of the graph model). 
This allows for graphs with non-negligible clustering and a specified degree distribution. 
 That is to say,  the CMC deviates from the classical  Erd\H{o}s-R{\'e}nyi graph model \citep{ER} in two fundamental ways: it allows for  
 for a non-Poissonian degree distributions and is asymptoticly clustered. 
 Epidemics on this type of graph have previously been studied by
 \citet{model_miller} and \citet{ODEs}.
\citet{model_miller} investigated the impact of clustering on the epidemic threshold, formulated as a bond percolation problem.  
This means that the infectivity of infected individuals is assumed to be homogeneous; an infected individual transmits the disease to each of its neighbors independently with some fixed probability $T$. 
\citet{ODEs} investigated the time evolution and final size of epidemics on CMC graphs  under the assumption 
of exponentially distributed infectious periods during which
individuals contact neighbors at a constant rate.

The main contribution of our research is that we extend the results of \citet{model_miller} and \citet{ODEs} by allowing  
for heterogeneous infectivity, i.e. by allowing for 
some infected individuals to be more contagious than others.
 Such heterogeneity may, for instance, reflect variability in the infectious period. 
 We provide expressions for the probability of a major outbreak and the final size of an major outbreak. 
A key tool in our analysis is the approximation of the epidemic seen from a ``generation of infection'' or ``rank'' perspective by a multitype Galton Watson branching process. 
This approximation, which is interesting in its own right, gives rise to the rank based reproduction number $R_0$ (see e.g.\ \citet{rank,R0I}). 

The second contribution of this paper concerns vaccination. We investigate the impact of uniform vaccination (i.e. vaccinated individuals are selected uniformly at random) with a perfect vaccine (i.e. a vaccine that provides full and permanent immunity to the disease). 
We find that it is necessary to vaccinate a fraction $1-1/R_0$ of the population  in order to prevent a major outbreak of the disease, as in the case of homogeneous mixing.
We illustrate our findings with numerical examples.

This paper is structured as follows. In Section \ref{sec:prel} we provide the preliminaries for the model. In Section \ref{sec:model} we give a more detailed description of how graphs are generated in the CMC and investigate the asymptotic clustering of such graphs and 
in Section \ref{sec:ep_mod} the epidemic model is specified. Section \ref{sec:R0}- \ref{sec:bp_approx} contains an overview of the concept of reproduction numbers and the necessary branching process background.  
In Section \ref{sec:not_vacc}, we derive expressions for 
the probability of a major outbreak and the expected final size
under the assumption of an unvaccinated and fully susceptible population, and in Section \ref{sec:vaccination} the analysis is repeated under the assumption of uniform vaccination with a perfect vaccine. We illustrate our findings with numerical examples presented in Section \ref{sec:num} and discuss possible extensions in Section \ref{sec:discussion}. 

\section{Preliminaries}\label{sec:prel}
\subsection{The configuration model with clustering}\label{sec:model}

\begin{figure}[h]
\centering

\includegraphics[]{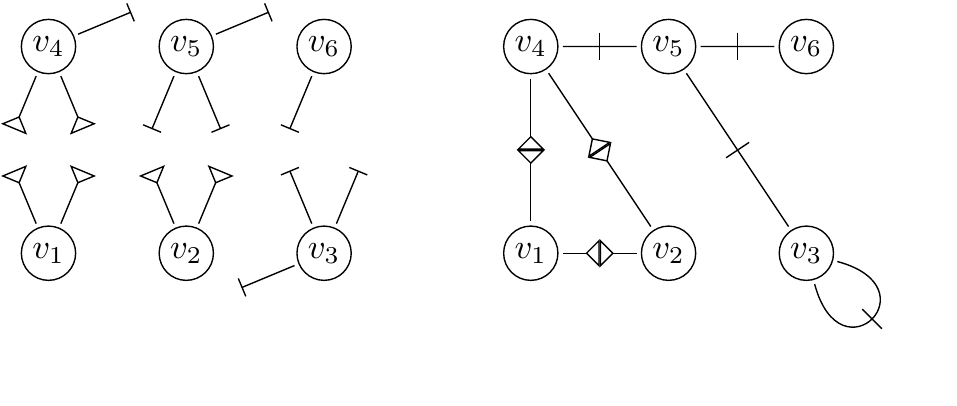}

\caption{Schematic illustration of the construction of a CMC graph. Triangle half-edges (marked with a triangle) and single half-edges (marked with a perpendicular line) are assigned to the nodes of the graph (left). The half-edges are then matched uniformly at random (right). Note that two of the half-edges attached to $v_3$ are paired with each other and so form a self-loop.  } \label{fig:conf_mod}
\end{figure}

A  configuration model with clustering CMC graph is constructed as follows.
Let $\{p(k_s, k_\Delta)\}_{k_s,k_\Delta\in\mathbb{N}_0}$
be a prescribed joint degree distribution, where $k_s$ denotes the number of single edges attached to a node, and $k_\Delta$ denotes the number of pairs of triangle edges. Throughout, $(S,\Delta)$ is assumed to be a generic random vector  distributed according to $p$. Let $\{(S_i, \Delta_i)\}_{i=1}^N$ be a sequence of independent copies of $(S,\Delta)$. 
Analogously to the CM, a graph $G_N=G_N(p)$ of size $N$ is constructed by first assigning the single degree $S_i$ and the triangle degree $\Delta_i$ to the node $v_i$, $i=1,2,\ldots,N$. 
One may think of this step in terms of half-edges;  
to each node $v_i$, we attach $S_i$ single half-edges and $\Delta_i$ pairs of triangle half-edges.  
The single half-edges are then matched in pairs and the triangle half-edge pairs in threes by choosing a matching uniformly at random among all possible such matchings. 
The process of joining half-edges is illustrated in Figure \ref{fig:conf_mod}. 
As described in \citet{model_miller}, the matching may be carried out as follows.
Two lists of nodes, one single degree list and one triangle degree list are created. A node with joint degree $(k_s,k_\Delta)$ appears $k_s$ times in the single list and $k_\Delta$ times in the triangle list. The lists are then shuffled uniformly, and the nodes on positions $2m+1$ and $2m+2$ in the single degree list and positions $3m+1,\ 3m+2$ and $3m+3$ in the triangle degree list are matched, $m\in\mathbb{N}_0$.   

We define the \emph{total single degree}\label{tot_sin} as
$$D^{(N)}_S:=\sum_{i=1}^NS_i$$
and the \emph{total triangle degree}\label{tot_tri} as
$$D^{(N)}_\Delta:=\sum_{i=1}^N\Delta_i.$$

If the total single degree
(that is, the length of the single degree list) is not even or if the total triangle edge degree (the length of the triangle degree list) is not a multiple of three we erase a single half-edge and/or one or two triangle half-edge pairs chosen uniformly at random. Similarly, we erase self-loops and merge multiple edges, so that the resulting graph is simple. 
Under assumption \ref{AssumDp1} (stated below) on $p$ it holds that the number of single self-loops and single double edges converge in distribution to independent Poisson random variables with finite means (cf. \citet[Prop. 7.13]{remco}). 

For this reason, self-loops and multiple edges are negligible in the limit as $N\to\infty$. In the remainder of this paper, we ignore the small differences in the topology of the graph that arise from erasing multiple edges or self-loops. In addition, we ignore the small differences in effective degree distribution that arise from erasing half-edges so that the number of single and triangle half-edges are multiples of two and three, respectively.

We make the following assumptions on $p$.
\begin{enumerate}[label=A\arabic*), ref=A{\arabic*}]
\item $E(\Delta^2)<\infty\text{ and } E(S^2)<\infty. $
\label{AssumDp1}
\item $P(\max(\Delta,S)\geq 2)>0$ and $E( \Delta S)>0$.\label{AssumDp3}
\end{enumerate} 

Note that the assumption \ref{AssumDp1} implies  $E(\Delta S)<\infty$. 
Assumption \ref{AssumDp3} ensures that the mean matrices of the approximating branching processes (presented below) are positively regular (we say that an $r\times r$ matrix $M$ is positively regular if it has finite non-negative entries and for some $n\in\mathbb{N}$ all entries of $M^n$ are strictly positive). 

\subsubsection{Clustering coefficient of $G_N$}

For any undirected graph we can measure the ammount of clustering in the network using the so-called clustering coefficient, which is defined as follows.
Let $G=(V, E)$ be an undirected graph with node set $V$ and edge set $E$. Define
$$\mathcal W^G_\wedge=\{(u,v,w)\in V^3: (u,v), (v,w)\in  E\}$$
 the set of all ordered wedges (i.e. directed paths consisting of precisely two edges) of $G$ and 
$$\mathcal W_\Delta^G=\{(u,v,w)\in V^3: (u,v), (v,w), (w, u)\in E\}\subset \mathcal W_\wedge^G$$
the set of all ordered triangles of $G.$ 
The \textit{clustering coefficient} $C(G)$ of $G$ is a measure of the degree of clustering of $G$ and is defined as the fraction of the ordered wedges of $G$ that are also triangles:
$$ C(G)=\frac{\vert \mathcal W_\Delta^G\vert } {\vert \mathcal W^G_\wedge\vert }.$$

As stated in the following proposition, 
CMC graphs 
have asymptotically non-zero clustering
as $N\to\infty$. An analogous result for fixed degree sequences was presented in \citet{newmanconf}. Let $\overset{P}{\longrightarrow}$ denote convergence in probability.

\begin{proposition}\label{clust_asympt}
Let $\{G_N\}_N$ be a sequence of CMC graphs with independent degrees drawn from $p$. If $p$ satisfies assumption \ref{AssumDp1} then
\begin{gather}
C(G_N)\overset{P}{\longrightarrow} \frac{E(2\Delta)}{E((2\Delta+S)^2)-E(2\Delta+S)}. 
\end{gather} 
\end{proposition}

The proof is presented in the Appendix. 

\subsubsection{Downshifted size-biased degrees}






The graph $G_N$ may be constructed by joining the half-edges in a random order. In particular, $G_N$ may be constructed as the epidemic progresses; starting with the initial infected case  we sequentially match the half-edges along which the disease is transmitted.  Since half-edges are chosen uniformly at random in the matching procedure, the probability to choose a specific node  
is proportional to the number of free half-edges attached to the node in question. That is, if we pair a single half-edge, the probability of choosing a specific node with $k_s$ unpaired single half-edges is proportional to $k_s$. 
For this reason, the degree distribution a node explored by joining a single half-edge in the early phase of the epidemic can be approximated by the \emph{single size biased} degree distribution $p^{(s)}$
\begin{gather}\label{s_size_biased_emp}
p^{(s)}(k_s, k_\Delta)=\frac{k_sp(k_s, k_\Delta)}{E(S)}.
\end{gather}

Similarly,  the degree distribution of the nodes explored by joining three triangle half-edge pairs in the early phase of the epidemic can be approximated by the \emph{triangle size biased} degree distribution $p^{(\Delta)}$
\begin{gather}\label{tri_size_biased_emp}
p^{(\Delta)}( k_s, k_\Delta)=\frac{k_\Delta p(k_s, k_\Delta)}{E(\Delta)}.
\end{gather}

In the epidemic process, we need to account for the fact that an infected individual has at least one non-susceptible neighbor (namely the direct source of its infection). For this reason, 
we introduce the 
\emph{downshifted} size biased degree distributions $p^{(s)}_\bullet$ and $ p^{(\Delta)}_\bullet$,
given by
\begin{gather}\label{s_size_biased_emp_shift}
\begin{aligned}
p^{(s)}_\bullet(k_s, k_\Delta)&=p^{(s)}( k_s+1, k_\Delta)\\
p^{(\Delta)}_\bullet(k_s, k_\Delta)&=p^{(\Delta)}(k_s, k_\Delta+1).
\end{aligned}
\end{gather}

Throughout, we will make frequent reference to the following random vectors
\begingroup
\addtolength{\jot}{.5em}
\begin{gather}\label{rvs}
\begin{aligned}[c]
(S^{(s)}_\bullet, \Delta^{(s)}_\bullet)&\sim p_\bullet^{(s)}\\
(S^{( \Delta)}_\bullet, \Delta^{(\Delta)}_\bullet)&\sim p_\bullet^{(\Delta)}\\
\end{aligned}
\end{gather}
\endgroup
and the expected values
\begingroup
\addtolength{\jot}{.5em}
\begin{gather}\label{rvs_exp}
\begin{aligned}[c]
E(S^{(s)}_\bullet)& =\frac{E(S^2)}{E(S)}- 1\\
E(S^{( \Delta)}_\bullet)&=\frac{E(S\Delta)}{E(\Delta)}.\\
\end{aligned}
\qquad
\begin{aligned}[c]
E(\Delta^{(s)}_\bullet)&=\frac{E(S\Delta)}{E(S)}\\
E(\Delta^{(\Delta)}_\bullet)&=\frac{E(\Delta^2)}{E(\Delta)}- 1\\
\end{aligned}
\end{gather}
\endgroup

\subsection{The epidemic model}\label{sec:ep_mod}

We use an SIR model to investigate the dynamics of the spread of the disease.  At any given time point, the population is divided into three groups, depending on health status. The groups are susceptible (\textbf{S}), infectious (\textbf{I}) and recovered (\textbf{R}) (see e.g.\ \citet{survey}).  
Individuals of the population make contact with other individuals at (possibly random) points in time.  
If, at some time point, an infectious individual contacts a susceptible individual then the susceptible individual instantaneously becomes infectious. An infectious individual will cease to be contagious after a period of time, which we call the \emph{infectious period} of the individual in question, and is then transferred to the recovered group. 
Recovered individuals are those that are immune to the disease. Individuals belonging to this group play no further role in the spread of the disease.
Because of this last observation, we can treat individuals that die because of the disease as ``recovered''.
In summary, we allow only the transitions $S\to I$ and $I\to R$. Note that the population is assumed to be closed; we ignore births, deaths 
 and migration.

More specifically, we consider an SIR epidemic in a generation framework on the clustered graph $G_N$ and assume heterogeneity in infectivity. 
That is, some infected individuals are more contagious than others. Such heterogeneity may, for instance,  
arise from variability in the infectious period. 
To this end, let $T$ be a random variable with support in $ [0,1]$, and let $\{T_i\}_{i=1}^N$ be a sequence of independent  copies of $T$.
Each node $v_i$ of $G_N$ is equipped with a \emph{transmission weight} $T_i$.  
If $v_i$ gets infected, then each susceptible neighbor of $v_i$ gets infected by $v_i$ independently  in the next generation with probability $T_i$ (conditioned on $\{T_i\}_i$). Node $v_i$ thereafter becomes recovered, playing no further role in the epidemic. 
An infected node transmits the disease independently of the transmissions from other infected nodes.  
An infected node does not, however, transmit the disease to its neighbors independently,
unless the distribution of $T$ is degenerate. 
Conditioned on the transmission weights $\{T_i\}_i$ and the structure of $G_N$, the number of neighbors that an infected node $v_i$ makes (potentially infectious) contact with while infectious has a binomial distribution with parameters $d_i$ and $T_i$, where $d_i$ is the degree of $v_i$.

The spread of this epidemic can be fully captured by a directed graph (see e.g.\ \citep{R0I}). To construct such directed graph from an undirected CMC graph $G_N$, we replace each undirected edge of $G_N$ by two parallel directed edges, pointing in the opposite direction.
 The weight of an edge $(v_i,v_j)$, which represents the (potential) transmission time from $v_i$ to $v_j$, is taken to be 1 if $v_i$ would make infectious contact with $v_j$ if infected, and $\infty$ otherwise. The individuals ultimately infected are then the individuals that can be reached from an initial case by following a path consisting of directed edges with finite edge weights.

\subsection{Reproduction numbers}\label{sec:R0}

A key quantity in the study of epidemics is the basic reproduction number,
often denoted by $R_0$. It is 
usually defined as the expected number of infected cases caused by a ``typical'' infected individual in an otherwise susceptible population. 
 For most stochastic epidemic models (including SIR epidemics in homogeneous mixing propulations \citep{survey}, populations with households \citep{R0II} and epidemics on networks \citep{CM_vacc}) it has the threshold property that a major outbreak is possible if and only if $R_0>1$. 
 For models where a suitable generation based branching process approximation is available, $R_0$ is usually defined as the Perron root (the dominant eigenvalue, which exists and is real-valued by assumptions \ref{AssumDp1} and \ref{AssumDp3}, see for instance \citet[Chapter 2]{varga}) of the mean matrix of the approximating Galton Watson branching process. 
This is the definition used in this article. 
By standard branching process theory, the interpretation of $R_0$ as the expected number
of cases caused by the typical individual in the early phase of the epidemic and its threshold properties are retained
by this definition. The threshold property of $R_0$ is made precise in Theorem \ref{major} below. 

In Section \ref{sec:vaccination}, we investigate the spread of an epidemic in a population with vaccination. 
To this end, in addition to the basic reproduction number $R_0$, we consider the \emph{perfect vaccine-associated reproduction number} $R_V$. 
A vaccine is \emph{perfect} if 
it provides full and permanent immunity. That is, an individual vaccinated with a perfect vaccine cannot contract the disease. 
The {perfect vaccine-associated reproduction number} $R_V$ is defined as \citep{R0II} 
\begin{gather}\label{RV}R_V=\frac{1}{1-f_{\text v}^{(c)}},
\end{gather}
where the \emph{critical vaccination coverage} $f_{\text{v}}^{(c)}$\label{crit_vacc_cov}  is the
fraction of the population that has to be vaccinated with a perfect vaccine in order to reduce $R_0$ to unity, if the vaccinated individuals are chosen uniformly at random. That is to say, $f_{\text{v}}^{(c)}=1-1/R_V$ is the fraction necessary to vaccinate in order to be guaranteed to prevent a major outbreak
\citep{survey}. Note that if $R_0\leq 1$ then $f_{\text v}^{(c)}=0.$

For many models, including 
epidemics on graphs generated by the CM \citep{CM_vacc}
and the standard stochastic SIR epidemic model (i.e. individuals mix homogeneously, see for instance \citet{survey}), 
$R_V=R_0$. That is,  
vaccinating a fraction $1-1/R_0$ of the population with a perfect vaccine is sufficient to surely prevent a major outbreak. 
On the other hand, for the households and households-workplaces model with uniform vaccination,
$R_V\geq R_0$ 
\citep{R0II} with strict inequality possible.
In Section \ref{vacc} we show that for the model analyzed in this report, 
$R_V=R_0.$

\subsubsection{Epidemics in continuous time - the rank based approach} \label{rank_based}

As mentioned above, heterogeneity in infectivity might arise from  heterogeneity in the infectious period; an important special case of the above described model is  epidemics in continuous time with random infectious periods where contacts between individuals take place according to point processes on $\mathbb{R}_{\geq 0}$.
Ignoring the real time-dynamics of an epidemic does not impact results that concern the final outcome of the epidemic.
This result was first presented by \citet{Ludwig}, see also \citet{rank} for a more recent discussion. This leads us to the often more tractable \emph{rank based} approach.

\begin{figure}[h]
\centering
\includegraphics[]{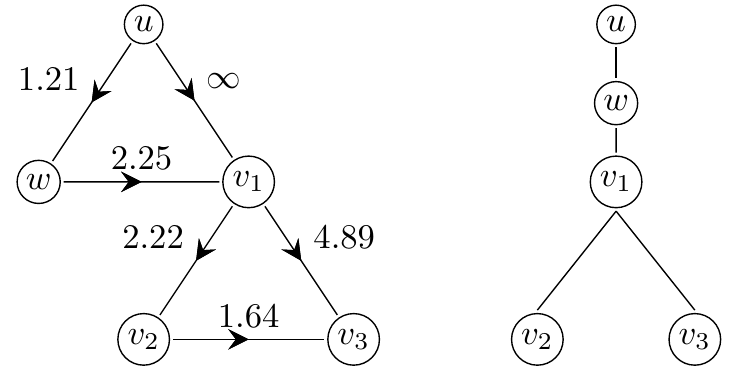}
\caption{{The difference between rank based generations and true generations}. Left: The length of the path $v_1\to v_3$ (i.e. the transmission time from $v_1$ to $v_3$) is $4.89$ and exceeds the length $2.22+1.64$ of the path $v_1\to v_2\to v_3.$ Therefore, the true path of transmission is $u\to w\to v_1\to v_2\to v_3.$ In the rank based approach, however, $v_3$ is attributed to $v_1$. Right: The resulting rank generation tree.   
} \label{fig:M3}
\end{figure}

In order to define the rank of a vertex, denote the initial case by $v_*$. The $rank$ of a node $v$ in $G_N$ is the distance from $v_ *$ to $v$, if every edge along which the disease would be transmitted is assigned the edge weight 1, and every other edge is assigned the edge weight $\infty$. That is, the rank of $v$ is the smallest number of directed edges that have to be traversed in order to follow a path of (potential) transmission from $v_*$ to $v$. 
We may then analyze the spread of the disease by letting generation $n$ of the epidemic process consist of the individuals of rank $n$. 
If, for instance, $v_1$ is the first node in a triangle consisting of the nodes $v_1,v_2,v_3$ to be infected, and $v_1$ infects $v_2$ and thereafter attempts to infect $v_3$, then $v_3$ is attributed to $v_1$ regardless of whether $v_1$ or $v_2$ infected $v_3$. This is illustrated in Figure \ref{fig:M3}. 

Consider a continuous time epidemic formulated as follows.
Suppose that each infected individual remains infectious for a (random) period of time. The infectious periods are distributed as the random variable $\tau$,  $\tau\sim F$,  and independent (but identically distributed) for different nodes. 
Suppose further that a node makes contact with each neighbor independently at a Poisson rate $\beta$ while infected, and that susceptible individuals are fully susceptible, so that each infectious-susceptible contact results in transmission. Without loss of generality we may assume $\beta=1$, since we may rescale time (and $F$ accordingly). The transmission weight $T$ is then distributed as
$1-e^{-\tau} $, and 
$E(T)=1-\mathcal L(1)$
and
$E(T(1-T))= \mathcal L(1)-\mathcal L(2)$ where  $\mathcal L(z)=\int_{\mathbb{R}_+}e^{-z x}dF(x)$\label{lapl} is the Laplace transform of the infectious period.

\subsection{Branching process approximations}\label{sec:bp_approx}

To analyze the spread of the disease in the early stages of the epidemic, we employ a multi-type branching process approximation. The graph $G_N$ may be constructed by joining the half-edges in any suitable (possibly random) order. In particular, the graph $G_N$ may be constructed (or explored) as the epidemic progresses; starting with the initial infected case $u^*$ we sequentially match the half-edges along which the disease is transmitted. In the early phase of the epidemic, short cycles (except for the triangles formed by triangle edges) are unlikely to occur. 
For these reasons, the early spread of the disease is well approximated by a suitably chosen branching process.

Similarly, a branching process approximation 
can be used to approximate the expected final size of the epidemic 
\citep{backwardsimple,  backwardhousehold, backwardadv}. 
In the graph representation of an epidemic, an individual  contracts the disease if and only if there is a path of directed edges with finite edge weights from the initial case to the node representing the individual in question.

\begin{figure}[h]
\centering
\includegraphics[]{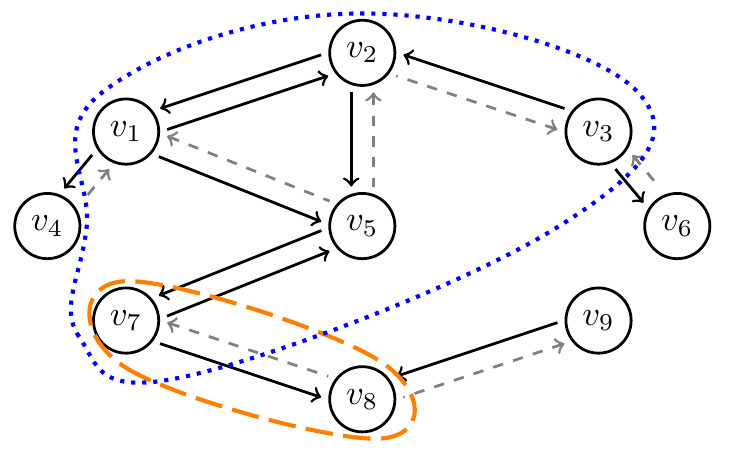}
\caption{Graph representation of an epidemic in a small $(N=9)$ population. The gray dashed and  black solid edges have infinite and finite edge weights (transmission times), respectively. 
The nodes in the susceptibility set  of $v_5$, $\mathcal S(v_5)=\{v_1.v_2,v_3,v_5,v_7\}$, are enclosed by the blue dotted line. The nodes that $v_5$ would infect if infected, directly or through other nodes, are enclosed by the orange dashed line.} \label{fig:susc}
\end{figure}

Define the susceptibility set  $\mathcal S(v)=\mathcal S_N(v)$ of a node $v$ as the collection of nodes of $G_N$ that can be reached from $v$ by tracing a path of finite length backwards.  
That is, the individuals that contract the disease are precisely the individuals 
with susceptibility sets that contain an initial case.
Hence, if the initial case is chosen uniformly at random then the probability that a node $v$ contracts the disease is proportional to the size of its susceptibility set $\mathcal S(v)$ and this probability can be approximated by exploring $\mathcal S(v)$. 
 Figure (\ref{fig:susc}) shows a schematic illustration of a susceptibility set.

By reversing the direction of the edges of the graph representation of an epidemic, but keeping the weights, the expected final fraction of the population infected in 
a major outbreak and the probability of a major outbreak are interchanged \citep{boundingmiller}, provided that the initial case is chosen uniformly at random. 
The process so obtained is called the \emph{backward epidemic process} 
of the node $v$.
If the underlying epidemic model is such 
that the {backward epidemic process} can be well approximated by a branching process, then we can use this branching process  
to compute the asymptotic distribution of the proportion 
of the population that ultimately escapes infection. 
This is made precise in the following theorem, due to \citet[Theorem 3.5]{backwardadv}, who proved the theorem for the related model of random intersection graphs.
The statement of Theorem \ref{major} carries over 
to the forward and backward branching processes considered in this paper. 
We omit the proof, which is analogous to the proof presented in \citet{ backwardadv}, see also \citet{backwardsimple}. Let $\overset{d}{\to}$ denote convergence in distribution.

\begin{theorem}\label{major}
Let $q$ and $q_b$ be the extinction probabilities of the forward and backward  approximating branching processes respectively, and let $S_N$ be the proportion of the population that ultimately escapes an epidemic in a population of size $N$. Then   
$$S_N\overset{d}{\to} S$$
	as $N\to\infty$ where $P(S = 1)=1-P(S = q_b)=q$. 
\end{theorem}

In other words, in the limit of large population sizes, the epidemic ``takes off'' with probability $ 1-q$, and if this happens a fraction $ 1-q_b$ of the population is ultimately infected (with probability converging to 1 as $N \to \infty$). Note that since $R_0$ is defined as the Perron root of the mean matrix of the forward branching process, $q<1$
 if and only if $R_0>1$.

\section{An epidemic in a fully susceptible population}\label{sec:not_vacc}

We now have the tools to analyze the spread of an infectious disease on a graph generated by the CMC. In the present section, the population is assumed to be fully susceptible to the disease, apart from the initially infectious individuals. 

\subsection{Forward process}

Before analyzing the forward process, we need to set  some terminology. For a given triangle $u,v,w$, where $u$ is the first individual to be infected in the triangle $u,v,w$, we refer to $v$ and $w$ as \emph{twins.}
We approximate the spread of the disease during the early phase 
by a multi-type branching process consisting of the following three types (except for the initial case):

\begin{enumerate}[]
\item[Type 1:]  A node infected along a triangle whose twin is infected at the same time step or earlier 
\item[Type 2:] A node infected along a triangle edge that is not of type 1
\item[Type 3:]  A node infected along a single edge
\end{enumerate}

\begin{figure}[h]
\centering
\includegraphics[]{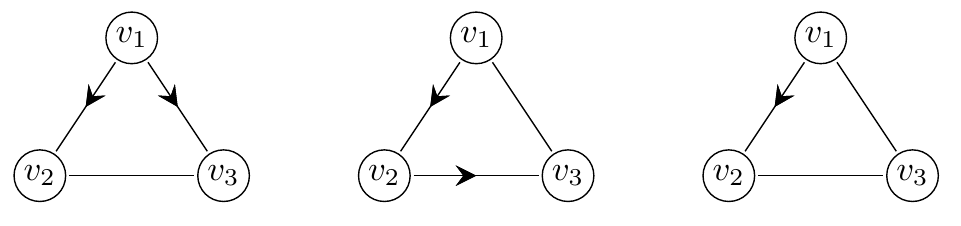}
\caption{Three examples of possible paths of transmission in a triangle $v_1,v_2,v_3$, where $v_1$ is the first node to be infected. Left: $v_1$ infects both $v_2$ and $v_3$. Both $v_2$ and $v_3$ are represented by type 1 individuals in the approximating branching process. Center: $v_1$ infects $v_2$ and $v_2$ infects  $v_3$. Then $v_3$ and $v_2$ are represented by type 1 and type 2 individuals, respectively. Right: $v_1$ infects $v_2$. Then $v_2$ is represented by a type 2 individual. } \label{fig:M4}
\end{figure}

Figure \ref{fig:M4} shows three examples of possible paths of transmission within a triangle giving rise to type 1 and 2 individuals in the approximating branching process.

Denote by 
$$M_f=(m_{ij})_{i,j=1}^3 = \begin{pmatrix}
m_{11} & m_{12} & m_{13}\\
m_{21} & m_{22} & m_{23}\\
m_{31} & m_{32} & m_{33}\\
\end{pmatrix}$$
the mean matrix of the above described branching process. 
Suppose that $v_1$ is the first individual to be infected in the triangle $v_1$, $v_2$, $v_3$. The probability that $v_1$ transmits the disease both to $v_2$ and $v_3$ is $E(T^2). $
Similarly, the probability that $v_1$ transmits the disease to either $v_2$ or $v_3$, but not to both, is $2E(T(1-T)). $

Thus, by linearity of expectation and because the distribution of the susceptible neighbors of infected nodes in the early phase of the epidemic is given by the downshifted degree distributions in (\ref{s_size_biased_emp_shift}), we obtain

\begin{equation}\label{mean_matr1}
\begingroup
\renewcommand*{\arraystretch}{2.5}
M_f=\begin{pmatrix}
2E(T^2)E(\Delta_\bullet^{(\Delta)}) & 2E(T(1-T))E(\Delta_\bullet^{(\Delta)})& E(T)E(S_\bullet^{(\Delta)})\\
2E(T^2)E(\Delta_\bullet^{(\Delta)})+E(T) & 2E(T(1-T))E(\Delta_\bullet^{(\Delta)}) & E(T)E(S_\bullet^{(\Delta)}) \\
2E(T^2)E(\Delta_\bullet^{(S)}) &   2E(T(1-T))E(\Delta_\bullet^{(S)}) & E(T)E(S_\bullet^{(S)})\\
\end{pmatrix}.
\endgroup
\end{equation}
(Recall that the random variables $\Delta_\bullet^{(\Delta)}$, $\Delta_\bullet^{(s)}$, $S_\bullet^{(\Delta)}$ and $S_\bullet^{(s)}$ defined in (\ref{rvs}) have the downshifted size biased distributions).
Note that all entries of $M_f$ are finite and that $S$ and $\Delta$ both have finite second moments by  assumption \ref{AssumDp1}.

If $M_f$ is positively regular (see the last paragraph in Section \ref{sec:model}) then $R_0$ is given by the Perron root of $M_f$. 
With little effort, one can use the expected values provided in (\ref{rvs_exp}) to show 
that necessary and sufficient conditions for $M_f$ to 
be positively regular are that assumptions \ref{AssumDp1}-\ref{AssumDp3} hold and that $0<E(T)<1$. 
If some of these conditions are not satisfied, we may analyze the spread of the disease by reducing the number of types of the approximating forward branching process. 

\subsubsection{Probability of a major outbreak}\label{sec:extinction}

For two $s$-dimensional vectors $\bar a=\transp{(a_1,\ldots, a_s)}$ and $\bar b=\transp{(b_1,\ldots, b_s)}$, we define \label{vecprod} $$\bar a^{\pow \bar b}:= a_1^{b_1}\cdot\ldots\cdot a_s^{b_s}.$$
Let
$f:[0,1]^3\to \mathbb{R}^3$
be the probability generating function of
the offspring distribution of the three types in the approximating branching process. That is, for $\bar{z}=\transp {(z_1,z_2,z_3)}\in[0,1]^3$  the $i$th component of $f(\bar{z})$ is given by 
\begin{gather}\label{eq0}
f(\bar{z})_i=E\left(\bar z^{\pow \bar \xi_i}\right)
\end{gather}
where $\bar \xi_{i}=(\xi_{i,1}, \xi_{i,2},\xi_{i,3})$ is distributed as the offspring of a type $i$ individual, $i=1,2,3.$ 

Similarly, let 
$f_*:[0,1]^3\to \mathbb{R}$ be the probability generating function of
the offspring distribution of the initial case. If $\bar \xi=\transp{(\xi_{*,1},\xi_{*,2},\xi_{*,3})}$ is distributed as the offspring of the initial case, then $f_*$ is
given by
$$f_*(\bar{z})=E\left(\bar z^{\pow\bar \xi}\right).$$

For $i=1,2,3$, let $(S^{(i)}, \Delta^{(i)})$ be 
 the joint degree of a type $i$ case with offspring $(\xi_{i,1}, \xi_{i,2},\xi_{i,3})$ and transmission weight $T$. That is, 
$$(S^{(1)}, \Delta^{(1)})\overset{d}{=}(S^{(2)}, \Delta^{(2)})\overset{d}{=}(S^{(\Delta)}, \Delta^{(\Delta)})$$
and 
$$(S^{(3)}, \Delta^{(3)})\overset{d}{=}(S^{(s)}, \Delta^{(s)}).$$
Here $\overset{d}{=}$ denotes equality in distribution.
By conditional independence we have
$$E(z_1^{\xi_{i,1}}z_2^{\xi_{i,2}}z_3^{\xi_{i,3}})=E\left(E(z_3^{\xi_{i,3}}\vert T,S^{(i)},\Delta^{(i)})E(z_1^{\xi_{i,1}}z_2^{\xi_{i,2}}\vert T,S^{(i)},\Delta^{(i)})\right).$$
Conditioned on the transmission weight $T$ and the single degree $S^{(1)}$, $\xi_{1,3}$ has a binomial distribution with parameters $S^{(1)}$ and $T$. 
Thus 
\begin{gather*}
\begin{split}
E(z_{3}^{\xi_{1,3}}\vert T,S^{(1)},\Delta^{(1)})=&\sum_{k_0+k_1=S^{(1)}}{{S^{(1)}}\choose {k_1} } (Tz_3)^{k_1}(1-T)^{k_0}\\
=&(Tz_3+1-T)^{S^{(1)}}.
\end{split}
\end{gather*}
Similarly 
\begin{gather*}
\begin{split}
\ & E(z_1^{\xi_{1,1}}z_2^{\xi_{1,2}}\vert T,S^{(1)},\Delta^{(1)})\\
&=\sum_{k_0+k_1+k_2=\Delta^{(1)}-1}{{\Delta^{(1)}-1}\choose {k_0,k_1,k_2} }(1-T)^{2k_0} (2(1-T)Tz_2)^{k_1}(Tz_1)^{2k_2}\\
&=((1-T)^2+2T(1-T)z_2+T^2z_1^2)^{\Delta^{(1)}-1}.
\end{split}
\end{gather*}
Thus 
\begin{gather}\label{eq1}
\begin{aligned}
&E(z_1^{\xi_{1,1}}z_2^{\xi_{1,2}}z_3^{\xi_{1,3}})\\
&=E((Tz_3+1-T)^{S_\bullet^{(\Delta)}}((1-T)^2+2T(1-T)z_2+T^2z_1^2)^{\Delta_\bullet^{(\Delta)}})
\end{aligned}
\end{gather}
where $(\Delta_\bullet^{(\Delta)}, S_\bullet^{(\Delta)})$ is independent of $T$. 

Since the conditional offspring distribution of a type 2 individual is identical to the offspring distribution of a type 1 individual except that a type 2 individual may give birth to one additional type 1 individual with probability $T$, we have
\begin{gather}
\label{eq2}
\begin{aligned}
&E(z_1^{\xi_{2,1}}z_2^{\xi_{2,2}}z_3^{\xi_{2,3}})\\
&=E((Tz_3+1-T)^{S_\bullet^{(\Delta)}}((1-T)^2+2T(1-T)z_2+T^2z_1^2)^{\Delta_\bullet^{(\Delta)}}(Tz_1+1-T)).
\end{aligned}
\end{gather}
Similarly, 
\begin{gather}\label{eq3}
\begin{split}
&E(z_1^{\xi_{3,1}}z_2^{\xi_{3,2}}z_3^{\xi_{3,3}})\\
&= E((Tz_3+1-T)^{S_\bullet^{(s)}}((1-T)^2+2T(1-T)z_2+T^2z_1^2)^{\Delta_\bullet^{(s)}}).
\end{split}
\end{gather}
Substituting (\ref{eq1})-(\ref{eq3}) into (\ref{eq0}) gives an expression for $f$.

 By standard branching process theory, if $R_0>0$  the extinction probability of a process descending from a type $i$ individual, $i=1,2,3,$ is given by 
$q_i$, where $\bar{q}=\transp {(q_1,q_2,q_3)}$ is the unique solution of 
$\bar{q}=f(\bar{q})$ in $[0,1)^3$. 
We also have
\begin{gather}\label{eq:test0}
\begin{split}
\bar q&=\lim_{n\to\infty} f^{\circ n}(\bar{0}), 
\end{split}
\end{gather}
where $f^{\circ n}$ is the composition of $f$ with itself $n$ times.

Since the approximating branching process dies out if and only if each of the processes started by the children of the initial case die out, the probability of extinction is given by 
$f^*(\bar{q}).$
After some calculations, analogous to the calculations that led to (\ref{eq1})-(\ref{eq3}), we find that the probability of extinction is given by 
$$f_*(\bar{q})=E\left((Tq_3+1-T)^{S}((1-T)^2+2T(1-T)q_2+T^2q_1^2)^{\Delta}\right)$$
where $(S, \Delta)$ is independent of $T$.  
We conclude that, by Theorem \ref{major}, the probability of a major outbreak is given by 
$1-f^*(\bar{q}),$
where $\bar q$ is the limit in (\ref{eq:test0}).

\subsection{Backward process}\label{sec_bp}

Let $w$ be a given node of $G_N$, chosen uniformly at random. We use a backward branching process to approximate the probability that $w$ contracts the disease, which by an exchangeability argument equals the expected final size of a major outbreak. 
The offspring of an individual $v$ in the backward process are the individuals that would potentially have infected $v$, if they were infected themselves.

The members of the susceptibility set are divided into the following two groups. This gives rise to a two-type approximating backward branching process.

\begin{enumerate}[]
\item[Type 1:] The vertex is included in the susceptibility set by virtue of potential transmission along a single edge
\item[Type 2:] The vertex is included in the susceptibility set by virtue of potential  transmission along a triangle edge
\end{enumerate}

We assign kinship as follows. The children of type 1 of an individual $v_1$ are the individuals included in the susceptibility set due to potential transmission along a single edge.
The children of type 2 of $v_1$ are the individuals included in the susceptibility set due to potential transmission of the disease to $v_1$, within a triangle of which $v_1$ is a member. 
We note that, given a triangle $v_1,v_2, v_3$ where $v_1$ is the primary case, both $v_2$ and $v_3$ will be members of the susceptibility set of $v_1$ by virtue of transmissions within the triangle if at least one of the following events happens:

\begin{enumerate}[label=$E_{\arabic*})$,ref=$E_{\arabic*}$]
\item $v_2$ and $v_3$ both ``infects'' $v_1$\label{e1}
\item $v_2$ infects $v_1$ and $v_3$ ``infects'' $v_2$\label{e2}
\item $v_3$ infects $v_1$ and $v_2$ ``infects'' $v_3$\label{e3}
\end{enumerate}
Here ``infects'' is conditional on the ``infector'' being infected during the epidemic.

\begin{figure}[h]
\centering
\includegraphics[]{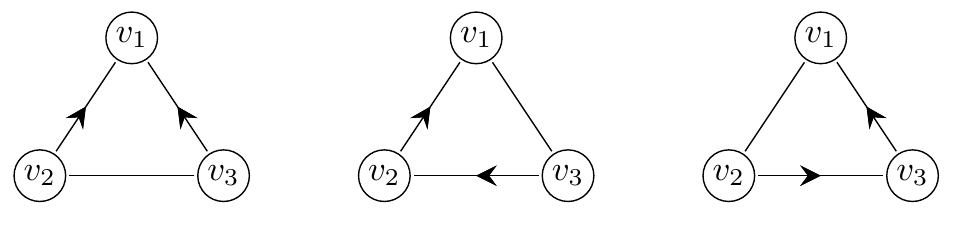}
\caption{The individuals $v_2$ and $v_3$ are both in the susceptibility set $\mathfrak{S}(v_1)$ of $v_1$ by virtue of transmission within the triangle $v_1,v_2,v_3$ if and only if at least one of the events \ref{e1} (left), \ref{e2} (center) or \ref{e3} (right) happens. } \label{fig:M5}
\end{figure}

The events \ref{e1}-\ref{e3} are illustrated in Figure (\ref{fig:M5}). 

Standard calculations give that the probability of the union of the events  \ref{e1}-\ref{e3} is given by $p_2=3E(T)^2-2E(T)E(T^2)$. Similarly, the probability that neither $v_1$ nor $v_2$ will be members of the susceptibility set of $v$ by transmissions within the triangle is given by $p_0=(1-E(T))^2.$ For later use, denote $1-p_0-p_2$ by $p_1$.

\subsubsection{Expected final size of a major outbreak}

Let $b$ be the probability generating function of the offspring distribution of the two types of the approximating backward branching process.
Furthermore, let $b_*$ be the probability generating function of
the offspring distribution of the ancestor $w$.
Analogously to the forward branching process, the probability that the bloodline started by a type $i,\ i=1,2$, individual will become extinct is given by $q^b_i$, where $\bar{q}_b=\transp {(q_1^b,q_2^b)}$ is the unique solution of $\bar{q}_b=b(\bar{q}_b)$ in $[0,1)^2$ (recall $R_0>1$). 
The probability of extinction is given by  $b_*(\bar{q}_b)$.

Proceeding in the same manner as in Section \ref{sec:extinction} yields
\begin{gather*}
b(z_1,z_2)_1=E\left((E(T)z_1+1-E(T))^{S_\bullet^{(s)}} (p_0+p_1z_2+p_2z_2^2)^{\Delta_\bullet^{(s)}}\right)
\end{gather*}
where 
$p_0,\ p_1$ and $p_2$ are as in Section \ref{sec_bp}. 
Similarly 
\begin{gather*}
b(z_1,z_2)_2=
E\left((E(T)z_1+1-E(T))^{S_\bullet^{(\Delta)}} (p_0+p_1z_2+p_2z_2^2)^{\Delta_\bullet^{(\Delta)}}\right),
\end{gather*}
and the probability of ultimate extinction of the backward process is given by 
\begin{gather*}
b_*(\bar{q}_ b)=E\left((E(T)q_1^b+1-E(T))^{S} (p_0+p_1q^b_2+p_2(q_2^b)^2)^{\Delta}\right).
\end{gather*}
We conclude that the expected final size of a major outbreak is given by 
$1-b_*(\bar{q}_ b). $

\section{Vaccination}\label{sec:vaccination}

\subsection{Random vaccination with a perfect vaccine}\label{vacc}

Assume that a fraction $f_{\text v}<1$ of the population is vaccinated, and that the vaccinated individuals are chosen uniformly at random (without replacement) from the population. The vaccine is perfect, in the sense that a vaccinated individual gains full and lasting immunity to the disease. 
If the population size $N$ is large, we may use a slightly different model, where
each individual is vaccinated with probability $f_{\text v}$, independently of the vaccination status of other individuals. By the law of large numbers, for our purposes the models are 
equivalent in the limit as the population size $N\to\infty. $

As before, we may approximate the early phase of the epidemic by a multi-type branching process.
The individuals of the approximating branching process are now of the following three types.

\begin{enumerate}[]
\item[Type 1:]  Infected along a triangle edge and has a twin that is known not to be susceptible
\item[Type 2:] Infected along a triangle edge and has a twin that might be susceptible
\item[Type 3:]  Infected along a single edge
\end{enumerate}

To clarify the types,
assume that in the early phase of the epidemic $v_1$ is the primary case in the triangle $v_1, v_2,v_3$. If $v_1$ attempts to transmit the disease both to $v_2$ and $v_3$ and succeeds (that is, none of $v_2$ and $v_3$ are vaccinated) then both $v_2$ and $v_3$ are represented by type 1 individuals in the approximating branching process. This happens with probability 
\begin{gather}\label{event11}
E(T^2)(1-f_{\text v})^2.
\end{gather}
If $v_1$ attempts to transmit the disease both to $v_2$ and $v_3$, but only succeeds to transmit the disease to $v_3$ (that is,  $v_2$ is vaccinated and $v_3$ is not vaccinated), then in the approximating branching process the individual representing $v_1$ gives birth to one type 1 individual (representing $v_3$) within the triangle $v_1,v_2, v_3$. This happens with probability 
\begin{gather}\label{event12}
E(T^2)f_{\text v}(1-f_{\text v}).\end{gather}
If $v_1$ attempts to transmit the disease only to $v_2$ and succeeds (that is,  $v_2$ is not vaccinated) then in the approximating branching process, the individual representing $v_1$ gives birth to one type 2 individual (representing $v_2$) within the triangle $v_1,v_2,v_3$. This happens with probability 
\begin{gather}\label{event13}
E(T(1-T))(1-f_{\text v}).\end{gather}

\begin{figure}[h]
\centering
\includegraphics[]{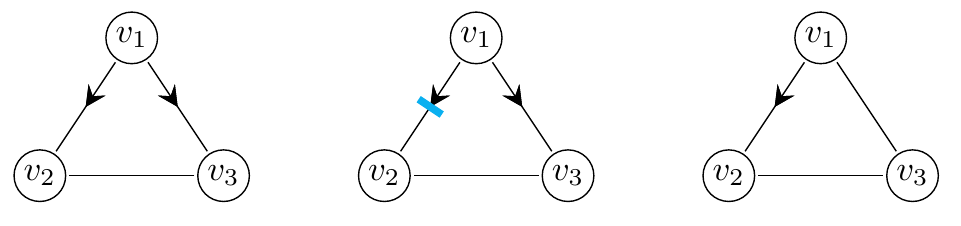}
\caption{Three examples of transmission dynamics within a triangle $v_1,v_2,v_3$.  An attempted transmission of the disease is represented by an arrow, an attempted transmission to a vaccinated individual is represented by an arrow and a blue bar. Left: $v_1$ attempts to transmit the disease both to $v_2$ and $v_3$, and succeeds. Both $v_2$ and $v_3$ are represented by type 1 individuals in the approximating branching process. Center:  $v_1$ attempts to transmit the disease both to $v_2$ and $v_3$, the transmission to $v_2$ is blocked since $v_2$ is vaccinated. Then $v_3$ is represented by a type 1 individual. Right: $v_1$ succeeds to transmit the disease to $v_2$, but does not attempt to infect $v_3$. Then $v_2$ is represented by a type 2 individual. } \label{fig:M8}
\end{figure}

The above described events are illustrated in Figure \ref{fig:M8}.

Denote the mean matrix of the approximating branching process by $M_ f^{(\text v)}=(m_{i,j}^{(\text v)})_{i,j=1}^3$. Using the expressions in (\ref{event11}) and (\ref{event12}) gives the expected number of type 1 individuals produced by a type 1 individual

\begin{gather*}
\begin{split}
m_{1,1}^{(\text v)}&=\left(2(1-f_{\text v})^2E(T^2)+2(1-f_{\text v})f_{\text v}E(T^2)\right)E\left(\Delta_\bullet^{(\Delta)}\right)\\
&=(1-f_{\text v})2E(T^2)E(\Delta_\bullet^{(\Delta)})\\
&=(1-f_{\text v})m_{1,1}
\end{split}
\end{gather*}

where $m_{1,1}$ is an element of the mean matrix $M_f$ of the forward branching process presented in (\ref{mean_matr1}).

Proceeding in the same fashion, we obtain the elements of the mean matrix $M_ f^{(\text v)}=(m_{i,j}^{({\text v})})_{i,j=1}^3$ of the branching process with random vaccination. It turns out that 
\begin{gather*}
\begin{split}
M_ f^{(\text v)}=(1-f_{\text v})M_f.
\end{split}
\end{gather*}
It is readily verified that the Perron root of $M_f^{(\text v)}$ is 
\begin{gather}\label{perron_eq}
r_f^{(\text v)}=(1-f_{\text v})r_f,
\end{gather}
where $r_f$ is the Perron root of $M_f$. Setting $r_f^{(\text v)}$ to 1 in (\ref{perron_eq}) and solving for $f_{\text v}$ yields the {critical vaccination coverage}
$
f_{\text v}^{(c)}=1-1/r_f.
$

We conclude that, for this particular graph model, equality holds between the basic reproduction number $R_0$ and the perfect vaccine-associated reproduction number $R_V$ as defined in (\ref{RV}).

\subsubsection{Probability of a major outbreak}

Let 
$h$
be the probability generating function of
the offspring distribution of the three types in our model including vaccination. As in Section \ref{sec:extinction}, we use the probability generating function to approximate the probability of extinction of the epidemic. 
To this end, let $(\zeta_{i,1}, \zeta_{i,2},\zeta_{i,3})$ be distributed as the offspring of a type $i$ individual with transmission weight  
$T$, $i=1,2,3, $ and let $(S^{(i)}, \Delta^{(i)})$ be distributed as the joint degree of this individual. 
That is, 
$$(S^{(1)}, \Delta^{(1)})\overset{d}{=}(S^{(2)}, \Delta^{(2)})\overset{d}{=}(S_\circ^{(\Delta)}, \Delta_\circ^{(\Delta)})$$
and 
$$(S^{(3)}, \Delta^{(3)})\overset{d}{=}(S_\circ^{(s)}, \Delta_\circ^{(s)}).$$
Note that $(S^{(i)}, \Delta^{(i)})$ and $T$ are independent.

By conditional independence

\begin{gather*}
\begin{split}E\left(z_1^{\zeta_{1,1}}z_2^{\zeta_{1,2}}z_3^{\zeta_{1,3}}\right)=E\left(E\left(z_3^{\zeta_{1,3}}\vert S^{(1)}, \Delta^{(1)}, T \right)E\left(z_1^{\zeta_{1,1}}z_2^{\zeta_{1,2}}\vert S^{(1)}, \Delta^{(1)}, T\right)\right)\end{split}
\end{gather*}
for $\bar{z}=\transp {(z_1,z_2,z_3)}\in[0,1]^3$.

Conditioned on the 
transmission weight $T$ and the joint degree $(S^{(1)}, \Delta^{(1)})$, the number of attempted transmissions from a type 1 individual along single edges has a binomial distribution with parameters $S^{(1)}$ and $T$, and each attempted transmission succeeds with probability $(1-f_{\text v})$. Thus, 

\begin{gather}\label{vac_eq11}
\begin{split}
E\left(z_3^{\zeta_{1,3}}\vert S^{(1)}, \Delta^{(1)}, T \right)&=\sum_{k_0+k_1=S^{(1)}} {S^{(1)} \choose  k_1}z_3^{k_1}\big(T(1-f_{\text v})\big)^{k_1}\big((1-T)+Tf_{\text v}\big)^{k_0}\\
&=\big(T(1-f_{\text v})z_3+1-T+Tf_{\text v}\big)^{S^{(1)}}.
\end{split}
\end{gather}

Similarly, for a type 1 individual $w$ with triangle degree $\Delta^{(1)}$, 
by conditioning on the number of attempted transmissions (in $k_i$ of the $\Delta^{(1)}-1$ triangles that are not yet affected by the disease, $w$ attempts to transmit the disease to $i$ individuals, $i=0,1,2$) and the vaccination status of the individuals  contacted by $w$ 
we obtain

\begin{gather}\label{vac_eq12}
\begin{split}
&E(z_1^{\zeta_{1,1}}z_2^{\zeta_{1,2}}\vert S^{(1)}, \Delta^{(1)}, T)\\[1em]
&=\sum_{k_0+k_1+k_2=\Delta^{(1)}-1}
{\Delta^{(1)}-1 \choose k_0, k_1, k_2}(1-T)^{2k_0}\big(2T(1-T)\big)^{k_1}T^{2k_2}\\[1em]
&\hphantom{{} \sum_{k_0+k_1+k_2=K_\Delta^1-1}} \p{\sum_{\tilde k_0+\tilde k_1+\tilde k_2=k_2} {k_2 \choose \tilde k_0, \tilde k_1,\tilde k_2} \big((1-f_{\text v})z_1\big)^{2\tilde k_2}\big(2f_{\text v}(1-f_{\text v})z_1\big)^{\tilde k_1}f_{\text v}^{2\tilde k_0}} \\[1em]
&\hphantom{{} \sum_{k_0+k_1+k_2=K_\Delta^1-1}} \p{\sum_{ k_0'+ k_1'=k_1} {k_1 \choose  k_0',  k_1'}(1-f_{\text v})^{k_1'}z_2^{k_1'} f_{\text v}^{k_0'}}\\[1em]
&=\sum_{k_0+k_1+k_2=\Delta^{(1)}-1}
{\Delta^{(1)}-1 \choose k_0, k_1, k_2}(1-T)^{2k_0}\big(2T(1-T)\big)^{k_1}T^{2k_2}\\[0.5em]
&\hphantom{{} \sum_{k_0+k_1+k_2=K_\Delta^1-1}} \big(\big((1-f_{\text v})z_1\big)^2+2f_{\text v}(1-f_{\text v})z_1+f_{\text v}^2\big)^{k_2} \\[0.5em]
&\hphantom{{} \sum_{k_0+k_1+k_2=K_\Delta^1-1}} \big((1-f_{\text v})z_2+f_{\text v}\big)^{k_1}\\[0.5em]
&=\Big[(1-T)^2 +2T(1-T)\big[(1-f_{\text v})z_2+f_{\text v}\big] \\
&+T^2\big[\big((1-f_{\text v})z_1)^2+2f_{\text v}(1-f_{\text v})z_1+f_{\text v}^2\big]\Big]^{\Delta^{(1)}-1}.\\
\end{split}
\end{gather}

Combining (\ref{vac_eq11}) and (\ref{vac_eq12}) yields

\begin{gather}\label{prob_gen1}
\begin{split}
E\left(z_1^{\zeta_{1,1}}z_2^{\zeta_{1,2}}z_3^{\zeta_{1,3}}\right)=E&\bigg[\Big(T(1-f_{\text v})z_3+1-T+Tf_{\text v}\Big)^{S_\bullet^{(\Delta)}}\\
&\hphantom{{} (} \Big((1-T)^2 +2T(1-T)\big((1-f_{\text v})z_2+f_{\text v}\big)\\
&\hphantom{{} (}\hphantom{{} ((1-T)^2 }+T^2\big(\big((1-f_{\text v})z_1\big)^2+2f_{\text v}(1-f_{\text v})z_1+f_{\text v}^2\big)\Big)^{\Delta_\bullet^{(\Delta)}}\bigg].
\end{split}
\end{gather}

By noting that the offspring distribution of a type 2 individual is identical to the offspring distribution of a type 1 individual, except that a type 2 may give birth to one additional type 1 individual with probability $T(1-f_{\text v})$ we obtain
\begin{gather}\label{prob_gen2}
\begin{split}
E\left(z_1^{\zeta_{2,1}}z_2^{\zeta_{2,2}}z_3^{\zeta_{2,3}}\right)=E\bigg[&\Big(T(1-f_{\text v})z_3+1-T+Tf_{\text v}\Big)^{S_\bullet^{(\Delta)}}\\[0.2em]
&\Big((1-T)^2 +2T(1-T)\big((1-f_{\text v})z_2+f_{\text v}\big)\\[0.2em]
&\hphantom{{}((1-T)^2}+T^2\big(\big((1-f_{\text v})z_1\big)^2+2f_{\text v}(1-f_{\text v})z_1+f_{\text v}^2\big)\Big)^{\Delta_\bullet^{(\Delta)}}\\[0.2em]
& \Big(z_1T(1-f_{\text v})+1-T(1-f_{\text v})\Big)\bigg].
\end{split}
\end{gather}

Similarly,  
\begin{gather}\label{prob_gen3}
\begin{split}
E\left(z_1^{\zeta_{3,1}}z_2^{\zeta_{3,2}}z_3^{\zeta_{3,3}}\right)=E\bigg[&\Big(T(1-f_{\text v})z_3+1-T+Tf_{\text v}\Big)^{S_\bullet^{(s)}}\\[0.2em]
&\Big((1-T)^2 +2T(1-T)\big((1-f_{\text v})z_2+f_{\text v}\big)\\[0.2em]
&\hphantom{{}(1-T)^2}+T^2\big(\big(\big((1-f_{\text v})z_1)^2+2f_{\text v}(1-f_{\text v})z_1+f_{\text v}^2\big)\Big)^{\Delta_\bullet^{(s)}}\bigg].\\
\end{split}
\end{gather}

Combining these results yields the probability generating function $h$ of the offspring distribution of a type 1, 2, 3 individual respectively. That is, $h(\bar z)_1$ is given by 
(\ref{prob_gen1}), $h(\bar z)_2$ is given by 
(\ref{prob_gen2}) and $h(\bar z)_3$ is given by 
(\ref{prob_gen3}). 

The probability generating function $h^*$ of the initial case is given by
\begin{gather}\label{vacc_eq2}
\begin{split}
h^*(\bar z)
=&E(z_1^{\zeta_{*,1}}z_2^{\zeta_{*,2}}z_3^{\zeta_{*,3}})\\[0.2em]
=&E\bigg[\Big(T(1-f_{\text v})z_3+1-T+Tf_{\text v}\Big)^{S}\\[0.2em]
&\hphantom{{} E(}\Big((1-T)^2 +2T(1-T)\big((1-f_{\text v})z_2+f_{\text v}\big)\\[0.2em]
&\hphantom{{}E(}\hphantom{{}((1-T)^2}+T^2\big(\big((1-f_{\text v})z_1\big)^2+2f_{\text v}(1-f_{\text v})z_1+f_{\text v}^2\big)\Big)^{\Delta}\bigg].\\
\end{split}
\end{gather}
for $\bar{z}=\transp {(z_1,z_2,z_3)}\in[0,1]^3$, where $(S, \Delta)$ is distributed as the joint degree of the initial case and independent of $T$. The probability of extinction of the approximating branching process is given by $h^*(\bar q^{\ (\text v)}),$
where $\bar q^{\ (\text v)}$ is given by the point in $[0,1]^3$ closest to the origin that satisfies 
$\bar q^{\ (\text v)}=h(\bar q^{\ (\text v)}).$
Thus, by Theorem \ref{major} the probability of a major outbreak is $1-h_*(\bar q^{\ (\text v)}).$

\subsubsection{The backward process}

We now turn our attention to the backward process and final size of an epidemic in a population where a fraction $f_{\text v}$ is vaccinated with a perfect vaccine.
To this end, we introduce the following three types, where individuals are classified by their vaccination status and the type of the edge along which they would transmit the disease if infected.

\begin{enumerate}[ ]
\item[Type 1:] Transmits along triangle edge, no information on vaccination status is available
\item[Type 2:] Transmits along triangle edge and is known not to be vaccinated since it is successfully infected by its twin 
\item[Type 3:] Transmits along single edge, no information on vaccination status 
is available
\end{enumerate}

\begin{figure}[h]
\centering
\includegraphics[]{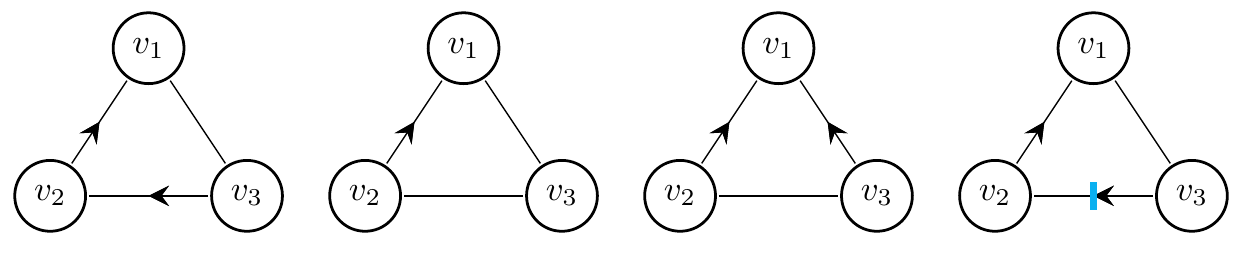}
\caption{At least one of $v_2$ and $v_3$ will belong to the susceptibility set of $v_1$ by virtue of potential transmissions within the triangle if some of the following types of scenarios (left to right in the picture) occur: \ref{Eb1}, \ref{Eb2}, \ref{Eb3}, \ref{Eb4}. An attempted transmission of the disease is represented by an arrow, an attempted transmission to a vaccinated individual is represented by an arrow and a blue bar.} \label{fig:M9}
\end{figure}

To clarify the types a bit more, let $v_1, v_2, v_3$ be a given triangle. 
At least one of $v_2$ and $v_3$ belongs to the susceptibility set of $v_1$ by virtue of potential transmissions within the triangle if some the following events, illustrated in Figure \ref{fig:M9}, happens. Note that all cases infected by virtue of transmission within the triangle $v_1,v_2, v_3$ are attributed to $v_1$.

\begin{enumerate}[label=$E_{\arabic*})$,ref=$E_{\arabic*}$]
\item $v_2$ attempts to infect $v_1$ and $v_3$ attempts to infect $v_2$, both succeed, and $v_3$ does not attempt to infect $v_1$. Or the same thing might happen, with $v_2$ and $v_3$ interchanged. This results in one type 1 and one type 2 individual in the approximating branching process. If $v_1$ is represented by a type 1 or 3 individual this happens with probability 
$$2\big(1-f_{\text v}\big)^2E(T)E\big(T(1-T)\big),$$
if $v_1$ is represented by an individual of type 2 this happens with probability
$$2(1-f_{\text v})E(T)E\big(T(1-T)\big).$$\label{Eb1}
\item  Only one of $v_2$ and $v_3$ attempts to infect $v_1$, and succeeds. The other node does not attempt to infect any node within the triangle. This results in one type 1 offspring. If $v_1$ is represented by an individual of type 1 or 3 this happens with probability 
$$2(1-f_{\text v})E(T)E\big(T(1-T)\big),$$
if $v_1$ is represented by an individual of type 2 this happens with probability
$$2E(T)E\big(T(1-T)\big).$$\label{Eb2}
\item  $v_2$ and $v_3$ both attempt to infect $v_1$ and succeeds. This results in two type 1 individuals born in the approximating branching process.  If $v_1$ is represented by an individual of type 1 or 3 this happens with probability $$(1-f_{\text v})E(T^2),$$
if $v_1$ is represented by an individual of type 2 this happens with probability
$$E(T^2).$$\label{Eb3}
\item $v_2$ attempts to infect $v_1$ and succeeds. The other node, $v_3$, attempts to infect $v_2$, but fails due to $v_2$ being vaccinated. The individual $v_3$ does not attempt to infect $v_1$. In this scenario, $v_2$ belongs to the susceptibility set of $v_1$. However, we do not include $v_2$ is the approximating branching process. This does not have any impact on the result of our analysis, since we are only interested in the probability of extinction of the backward process and $v_2$ does not produce any offspring in this process. \label{Eb4}
\end{enumerate}

\subsubsection{Expected final size}

Let $b^{(\text v)}$ and $b^{(\text v)}_*$ be the probability generating function of
the offspring distribution of the three types of the approximating backward branching process and of the ancestor, respectively. Furthermore, let $\bar \zeta_i=({\zeta^b_{i,1}},{\zeta^b_{i,2}}, {\zeta^b_{i,3}})$ be distributed as the offspring of a type $i,i=1,2, 3$, individual and denote by $E_s$ the conditional expectation given that the parent of $\zeta^b_{i,1},\zeta^b_{i,2}, \zeta^b_{i,3}$ is susceptible. 
 Let further $\bar \zeta_*=({\zeta^b_{*,1}},{\zeta^b_{*,2}}, {\zeta^b_{*,3}})$ be distributed as the offspring of the ancestor. Denote the extinction probability of a process descending from a type $i$ individual by $q_i^b$, $i=1,2,3$ and let $\bar{q}^b=\transp {(q_1^b,q_2^b, q_3^b)}$. 

To find an expression for $b^{(\text v)}$, we first note that for $\bar z=\transp{(z_1,z_2,z_3)}$
\begin{gather}\label{a1}
\begin{split}E\left(\bar z^{\pow \bar \zeta_3}\right)=f_{\text v}+(1-f_{\text v})E_s\left(E_s\left(z_3^{\zeta_{3,3}^b}\vert S^{(3)}, \Delta^{(3)} \right)E_s\left(z_1^{\zeta_{3,1}^b}z_2^{\zeta_{3,2}^b}\vert S^{(3)}, \Delta^{(3)}\right)\right)\end{split}
\end{gather}
where, as before, $(S^{(i)}, \Delta^{(i)})$ is distributed as the joint degree of a type $i$ individual, $i=1,2,3.$

Now
\begin{gather}\label{a2}
\begin{split}
E_s\left(z_3^{\zeta_{3,3}}\vert S^{(3)}, \Delta^{(3)}\right)&=\sum_{k_0+k_1=S^{(3)}-1} {S^{(3)}-1 \choose k_0, k_1}z_3^{k_1}E(T)^{k_1}E(1-T)^{k_0}\\
&=\big(E(T)z_3+1-E(T)\big)^{S^{(3)}-1}.
\end{split}
\end{gather}
By conditioning on the number of triangles $k_2$ in which an event of type \ref{Eb3} occurs, the number of triangles $k_1^a$ in which an event of type \ref{Eb1} occurs, the number of triangles $k_1^b$ in which an event of type \ref{Eb4} occurs and the number of triangles $k_1^c$ in which an event of type \ref{Eb2} occurs we obtain
\begin{equation}\label{a3}
\begin{split}
&E_s(z_1^{\zeta_{3,1}}z_2^{\zeta_{3,2}}\vert S^{(3)}, \Delta^{(3)})\\[0.5em]
&=\sum_{k_0+k_1^a+k_1^b+k_1^c+k_2=\Delta^{(3)}}{\Delta^{(3)} \choose k_0,k_1^a,k_1^b,k_1^c,k_2}E(1-T)^{2k_0}\\[0.5em]
&\hphantom{{}
=\sum_{k_0+k_1^a+k_1^b+k_1^c+k_2=K_\Delta^1}}\Big(2E(T)E\big(T(1-T)\big)(1-f_{\text v})\Big)^{k_1^a}\\[0.5em]
&\hphantom{{}
=\sum_{k_0+k_1^a+k_1^b+k_1^c+k_2=K_\Delta^1}} \Big(2E(T)E\big(T(1-T)\big)f_{\text v}\Big)^{k_1^b}\Big(2E(T)E\big((1-T)^2\big)\Big)^{k_1^c} \\[0.5em]
&\hphantom{{}
=\sum_{k_0+k_1^a+k_1^b+k_1^c+k_2=K_\Delta^1}} E(T)^{2k_2} z_2^{k_1^a}z_1^{k_1^a+k_1^c+2k_2}\\[0.3em]
&=\Big(\big(E(1-T)\big)^2+2E\big(T\big)E\big(T(1-T)\big)(1-f_{\text v})z_2z_1+2E(T)E\big(T(1-T)\big)f_{\text v}\\[0.3em]
&\hphantom{{} =((E(1-T))^2}+2E(T)E\big((1-T)^2\big)z_1+E(T)^2z_1^2\Big)^{\Delta^{(3)}}.
\end{split}
\end{equation}
Inserting the right hand sides of (\ref{a2}) and (\ref{a3}) in (\ref{a1}) gives
\begin{gather}\label{a4}
\begin{split}
E(z_1^{\zeta_{3,1}}z_2^{\zeta_{3,2}}z_3^{\zeta_{3,3}})\\[0.5em]
=f_{\text v}+(1-f_{\text v})E\bigg[&\Big(E(T)z_3+1-E(T)\Big)^{S_\bullet^{(s)}}\\[0.5em]
&\Big(\big(E(1-T)\big)^2 +2E(T)E\big(T(1-T)\big)(1-f_{\text v})z_1z_2\\[0.5em]
&\hphantom{((E(1-T))^2}+2E(T)E\big(T(1-T)\big)f_{\text v}\\[0.5em]
&\hphantom{((E(1-T))^2}+2E(T)E\big((1-T)^2\big)z_1+E(T)^2z_1^2\Big)^{\Delta_\bullet^{(s)}}\bigg].
\end{split}
\end{gather}
Similarly 
\begin{gather}\label{a5}
\begin{split}
E(z_1^{\zeta_{2,1}}z_2^{\zeta_{2,2}}z_3^{\zeta_{2,3}})=E\bigg[&\Big(E(T)z_3+1-E(T)\Big)^{S_\bullet^{(\Delta)}}\\[0.5em]
&\Big(\big(E(1-T)\big)^2 +2E(T)E\big(T(1-T)\big)(1-f_{\text v})z_1z_2\\[0.5em]
&\hphantom{((E(1-T))^2}+2E(T)E\big(T(1-T)\big)f_{\text v}\\[0.5em]
&\hphantom{((E(1-T))^2}+2E(T)E\big((1-T)^2\big)z_1+E(T)^2z_1^2\Big)^{\Delta_\bullet^{(\Delta)} }\bigg].
\end{split}
\end{gather}
and 
\begin{gather}\label{a50}
\begin{split}
E(z_1^{\zeta_{1,1}}z_2^{\zeta_{1,2}}z_3^{\zeta_{1,3}})=f_{\text v}+(1-f_{\text v})E(z_1^{\zeta_{2,1}}z_2^{\zeta_{2,2}}z_3^{\zeta_{2,3}}).
\end{split}
\end{gather}

Combining these results yields the probability generating function of the offspring distribution of the three types; $b^{(\text v)}(\bar z)_3$ is given by (\ref{a4}) and 
$b^{(\text v)}(\bar z)_2$ is given by (\ref{a5}). By replacing 
$(S^{(s)}_\bullet, \Delta^{(s)}_\bullet)$ in the right hand side of
(\ref{a4}) by $(S^{(\Delta)}_\bullet, \Delta^{(\Delta)}_\bullet)$ we obtain $b^{(\text v)}(\bar z)_1$.

Also by replacing $(S^{(s)}_\bullet, \Delta^{(s)}_\bullet)$ in the right hand side of
(\ref{a4}), but now by $(S,\Delta)$ we obtain the probability generating function $b_*^{(\text v)}(\bar z)$ of the offspring of the initial case. The expected final size of the epidemic, conditioned on that a major outbreak occurs, is given by $$1-b^{(\text v)}_*(\bar{q}^b).$$

\section{Numerical example}\label{sec:num}

Under very general assumptions, increasing the heterogeneity  in infectiousness leads to a decrease in the the probability of a major outbreak, the expected final size and $R_0$ \citep{ kuulasmaa,boundingmeester, boundingmiller}, see also \citet{ball_1985, secondlook, backward_miller}. 
In particular, for a fixed (marginal) transmission probability $E(T)$, the probability of a major outbreak and the expected final size are maximized  if $T=E(T)$ with probability 1 and minimized if $P(T=1)=E(T)=1-P(T=0)$. 
Similarly, for given $E(T)$, $R_0$ is maximized  if $T=E(T)$ with probability 1 and minimized if $P(T=1)=E(T)=1-P(T=0)$.

We illustrate this with the following example. 
Consider the three degree distributions
\begin{enumerate}[]
\item $p(2,1)=1$ \label{distr1}
\item $p(4,0)=0.95=1-p(2,1)$\label{distr2}
\item $p(0,2)=0.95=1-p(2,1)$\label{distr3}.
\end{enumerate}

That is, in all three degree distributions the total degree is 4 with probability 1. In addition,  distribution \ref{distr1} corresponds to a network where every node  is member of exactly one triangle. Distribution \ref{distr2} corresponds to a network where a node is not a member of any triangle with probability 0.95, while with probability 0.05 a node is member of one triangle. Finally, distribution \ref{distr3} corresponds to a network where a node is a member of two triangles with probability 0.95, while with probability 0.05 a node is member of one triangle.

\begin{figure}[h!]%
	\centering
	\subfloat[][]{%
		\label{fig:ex3-a}%
	\includegraphics[width=0.4\textwidth]{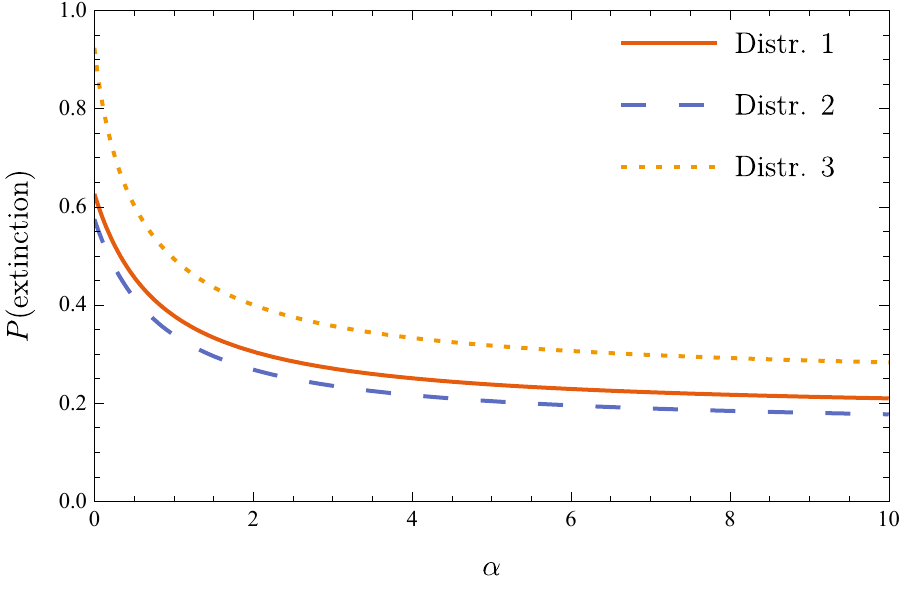}}%
	\hspace{8pt}%
	\subfloat[][]{%
		\label{fig:ex3-b}%
	\includegraphics[width=0.4\textwidth]{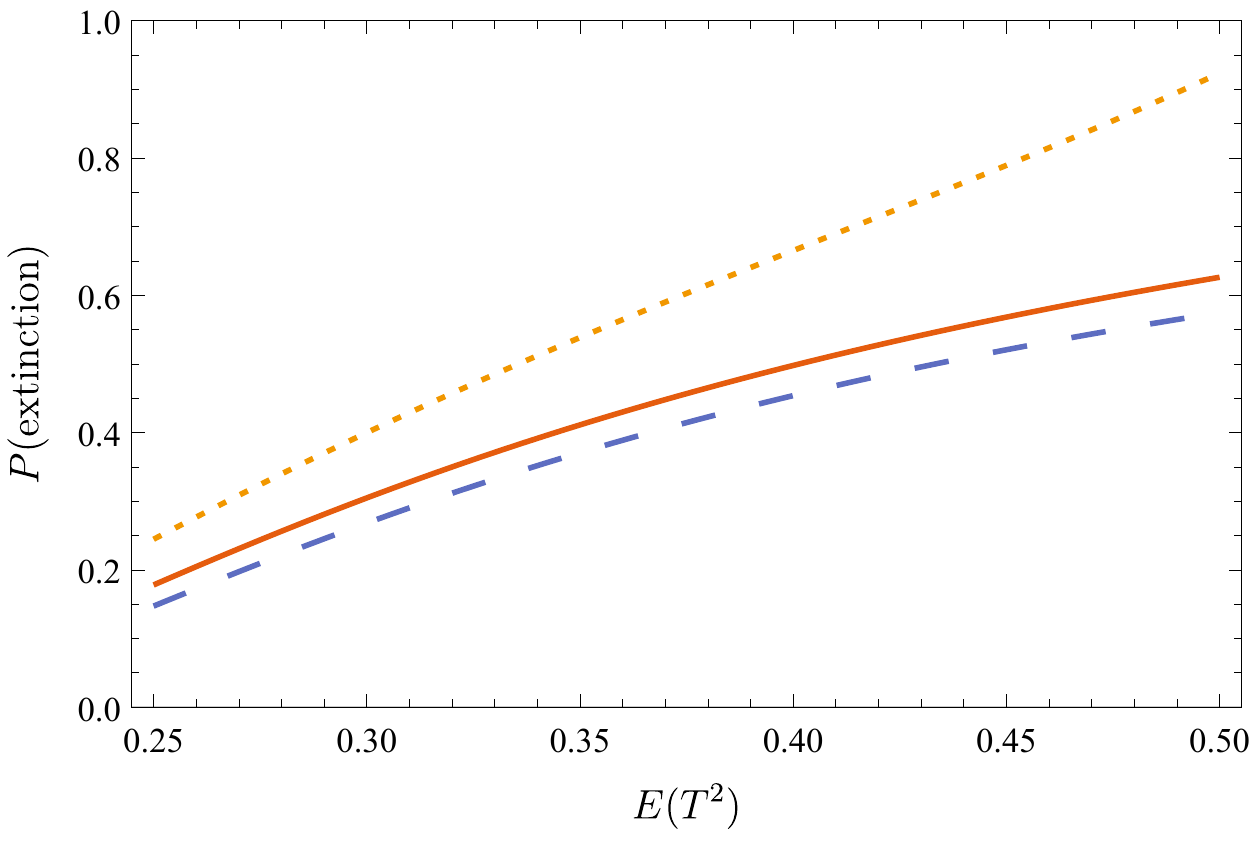}}\\
	\subfloat[][]{%
		\label{fig:ex3-c}%
	\includegraphics[width=0.4\textwidth]{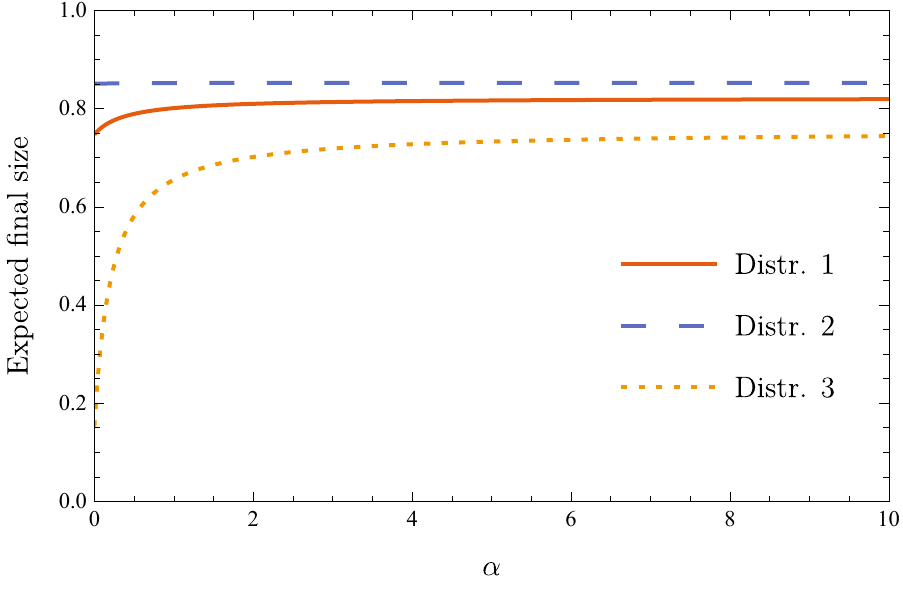}}%
	\hspace{8pt}%
	\subfloat[][]{%
		\label{fig:ex3-d}%
	\includegraphics[width=0.4\textwidth]{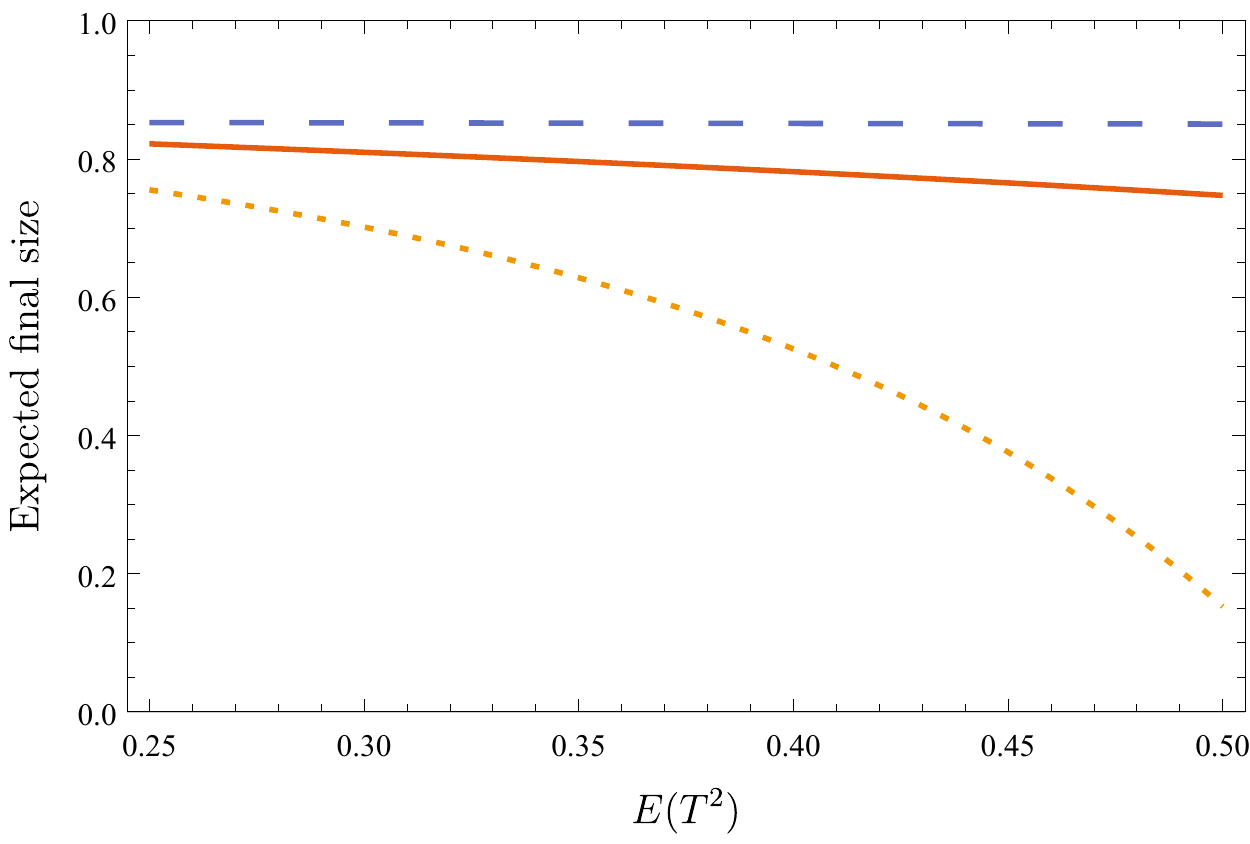}}\\
	\subfloat[][]{%
		\label{fig:ex3-e}%
		\includegraphics[width=0.4\textwidth]{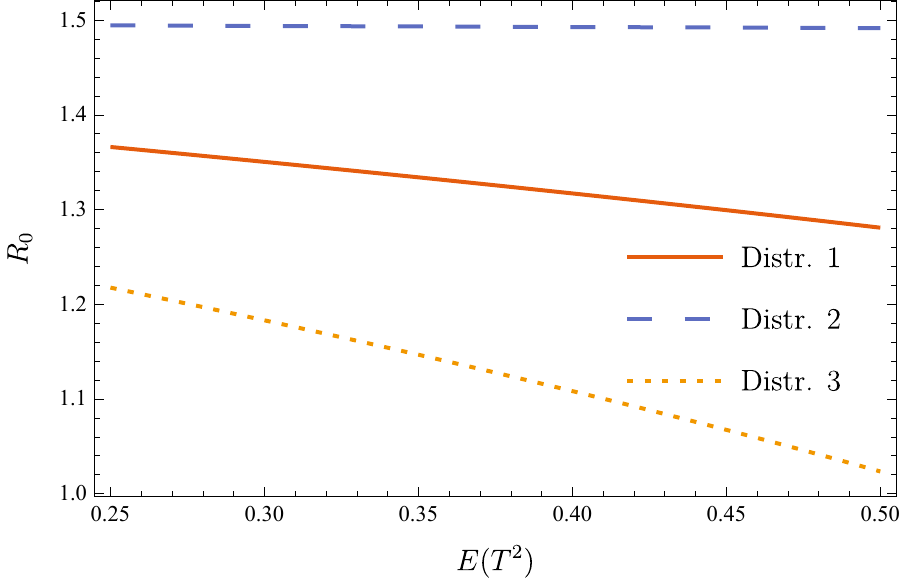}}%
	\hspace{8pt}%
	\subfloat[][]{%
		\label{fig:ex3-f}%
\includegraphics[width=0.4\textwidth]{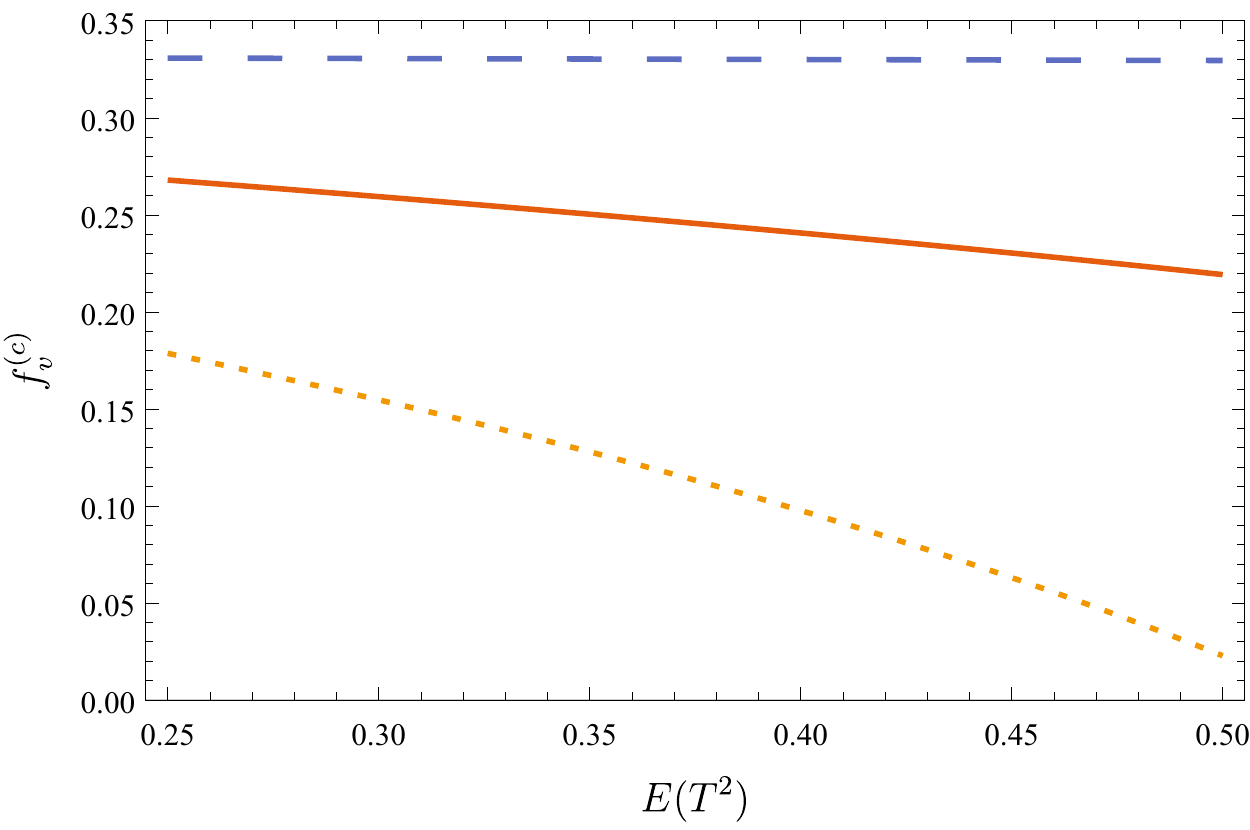}}\\
	\caption[something]{The impact of heterogeneity in infectivity for the three degree distributions.  
		\subref{fig:ex3-a} The probability that a major outbreak does not occur as a function of $\alpha$.
		\subref{fig:ex3-b} The probability that a major outbreak does not occur as a function of $E(T^2)$.
		\subref{fig:ex3-c} The expected final size of a major outbreak  as a function of $\alpha$.
		\subref{fig:ex3-d} The expected final size of a major outbreak  as a function of  $E(T^2)$.
		\subref{fig:ex3-e} The basic reproduction number $R_0$ as a function of $E(T^2)$.
		\subref{fig:ex3-f}The critical vaccination coverage $f_v^{(c)}$ as a function of $E(T^2)$.}%
	\label{fig3}%
\end{figure}

Furthermore, let $T$ have distribution $\text{Beta}(\alpha, \alpha)$ for some $\alpha>0$. That is, $T$ has density, $C_{\alpha}x^{\alpha-1}(1-x)^{\alpha-1}$, on the interval $(0,1)$, where $C_{\alpha}$ is a normalizing constant.  Then $E(T)=1/2$ and we can tune the heterogeneity of the infectivity of infected individuals by varying $\alpha$. In particular
$$E(T^2)=\frac{1}{2}\left(1-\frac{1}{2+\alpha^{-1}}\right).$$
Note that $\alpha\to\infty$ corresponds to $T$ being uniform on $(0,1)$, while $\alpha=0$ corresponds to $P(T=0)=P(T=1)=1/2.$
Figure \ref{fig3} shows
the 
probability that a major outbreak does not occur, the expected final size, $R_0$ and the critical vaccination coverage $f_v^{(c)}$
as functions of $\alpha$ or $E(T^2)$. 

As can be seen in Figure \ref{fig3}, ignoring actual heterogeneity  of infectivity in this case leads to an overestimation of the probability of a major outbreak (\ref{fig:ex3-a}-\ref{fig:ex3-b}). 
This effect is particularly evident in the presence of high clustering; the steeper slope of the curve corresponding to distribution 3 (\ref{fig:ex3-b}) 
and the relatively  low probability of a major outbreak when $\alpha$ is small
can be explained  by the fact that the approximating forward branching process is close to being critical when $\alpha$ is small. 
Figure \ref{fig:ex3-c}-\ref{fig:ex3-d} shows
that 
heterogeneity  of infectivity 
has virtually no impact on the expected final size of a major outbreak  and $R_0$ in the near absence of clustering.  
In the presence of clustering, on the other hand, ignoring heterogeneity of infectivity leads to an underestimation of the expected final size and a substantial overestimation of the critical vaccination coverage $f_v^{(c)}$. 
Note that $R_0$ and $f_v^{(c)}$ depend on the distribution of $T$ only through the first  and second moment of $T$.

\section{Discussion}\label{sec:discussion}


In this paper, we have incorporated clustering in the spread of an infectious disease by allowing for groups of size three with non-overlapping edges. It is, in principle, straightforward to extend the methods used in this paper to 
larger group sizes.
The CMC may, for instance, be generalized to larger group sizes as follows. %
Let $K=\{k_1,\ldots, k_r\}\subset\mathbb{N}_{\geq 2}$ be the set of possible group sizes. 
In the matching procedure, each node is equipped with an $r$-dimensional degree in $\mathbb{N}_0^r$. The $i$th component (the $k_i$-degree) of a degree specifies the number of groups of size $k_i$ to which the node in question belongs. 
Analogously to the construction of a CMC graph,
groups are then formed by creating one list for each group size; 
 a node with $k_i$-degree $d_i$ appears precisely $d_i$ times in the list corresponding to groups of size $k_i$. 
 The lists are then shuffled and half-edges of
nodes in positions $k+1,\ldots,k+k_i$ in the $k_i$-list are joined. 
The structure of a graph so obtained would be characterized by fully connected cliques, and similar to that of a random intersection graph \citep{backwardadv}. 
One possible approach to investigate epidemics on such graphs would be to approximate the spread of the disease by a multitype Galton Watson process where groups (or cliques) are represented by the particles of the branching process. The types of the approximating branching process   
would then be vectors in $\mathbb{N}^2$ of the form $(m,n)$, where $m$ represents the size of the clique and $n$ represents the number of members of the clique that the primary case of the clique attempts to infect. 
Another possible approach would be to use an infinite type 
branching process in the spirit of \citet{backwardadv}. 
We believe that the result would be analogous to the results obtained in \citet{backwardadv}. 

\appendix
\section*{Appendix: Proof of proposition \ref{clust_asympt}}\label{proof}

Let $\bar d=\{(S_i,\Delta_i)\}_{i\in\mathbb{N}}$ be a given (i.e. non-random) degree sequence that satisfies the following regularity  assumptions. 

\begin{enumerate}[label=A\arabic*), series=assumptions, ref=A{\arabic*}]
	\item $\frac{\sum_{i=1}^N\ind(S_i=k_1, \Delta_i=k_2) }{N }\to p(k_1,k_2)$ for any $k_1,k_2\in \mathbb{Z}_{\geq 0}$.	\label{R1}	
	\item $\frac{\sum_{i=1}^{N}\Delta_i^2}{N}\to E(\Delta^2) $ and  $\frac{\sum_{i=1}^{N}S_i^2}{N}\to E(S^2)$
	\label{R3}	
\end{enumerate} 

where $(S, \Delta)$ has distribution $p$, which is assumed to satisfy \ref{AssumDp1}-\ref{AssumDp3} in Section \ref{sec:model}. 
Let further $\overline G=\{G_N\}_{N\in\mathbb{N}}$ be a sequence of graphs generated by the CMC, where the degree sequence of $G_N$ is given by $\bar d_N=\{(S_i,\Delta_i)\}_{i=1}^N$ and denote $D_S^{(N)}=\sum_{i=1}^NS_i$.

Under the assumptions  \ref{R1}-\ref{R3} 
the expected number of self-loops and the expected number of multiple edges are borth of order $O(1)$
(cf. \citet[prop. 7.11]{remco}). 
Denote by $A_N$ the number of wedges of $G_N$ that are "deleted" when merging multiple edges and erasing self-loops, that is $$A_N=\sum_{i=1}^N{ S_i+2\Delta_i \choose 2}2-\vert \mathcal W^{G_N}_\wedge\vert=\sum_{i=1}^N (S_i+2\Delta_i)(S_i+2\Delta_i-1)-\vert \mathcal W^{G_N}_\wedge\vert,$$
then $E(A_N)=O(1)$.


From the definition of $A_N$, we deduce that
the total number of ordered triangles of $G_N$ is bounded from below by
$\vert \mathcal W_\Delta^{G_N}\vert\geq \sum_{i=1}^N  2\Delta_i-A_N$
and the total number of ordered wedges is bounded from above by
$$\vert \mathcal W^{G_N}_\wedge\vert\leq\sum_{i=1}^N{ S_i+2\Delta_i \choose 2}2=\sum_{i=1}^N (S_i+2\Delta_i)(S_i+2\Delta_i-1). $$
Therefore, by the definition of $C(G_N)$ and the assumptions above
\begin{gather}\label{clust_assympt}
C(G_N)\geq \frac{\p{\frac{\sum_{i=1}^N 2\Delta_i}{N} }-A_N} {\p{\frac{\sum_{i=1}^N (S_i+2\Delta_i)(S_i+2\Delta_i-1)}{N}}}\overset{P}{\to} \frac{E(2\Delta)}{E((2\Delta+S)^2)-E(2\Delta+S)}
\end{gather}
as $N\to\infty$.

This lower bound is tight in the limit as the number of nodes $N\to \infty$. 
Indeed, denote by $\mathcal W_s^{G_N}$ the set of ordered triangles of $G_N$ that consists solely of single edges, i.e. 
$$\mathcal W_s^{G_N}=\{(u,v,w)\in V_N^3:\ (u,v),(u,w)\text{ and $(v,w)$ are single edges} \},$$
where $V_N$ is the node set of $G_N$. 
Now, whenever $D_S^{(N)}\geq 6$
\begin{gather}\label{tri}
E \left(\vert \mathcal W_s^{G_N}\vert \right)\leq\sum_i\left( { S_i \choose 2} \sum_{j}\frac{S_j}{D_S^{(N)}-2}\left(\sum_{l}\frac{S_l}{D_S^{(N)}-3}\left(\frac{(S_j-1)(S_l-1)}{D_S^{(N)}-5}\right)\right)\right) 
\end{gather}
where the sums run over the integers $1,\ldots, N$.

Dividing by $N$ in (\ref{tri}) and letting $N$ approach infinity gives 
$E(\vert \mathcal W_S^{G_N}\vert)/N\to 0$
as $N\to\infty$.
Thus 
$\vert \mathcal W_S^{G_N}\vert/N{\to }\ 0$  
in probability. 
Repeating this procedure for triangles formed by a combination of triangle and single edges gives  
\begin{gather}\label{clust_lim}
C(G_N)\overset{P }{\longrightarrow} \frac{E(2\Delta)}{E((2\Delta+S)^2)-E(2\Delta+S)}. 
\end{gather} 
The assertion now follows by bounded convergence and the law of large numbers.

\noindent \subsection*{Acknowledgements}
We thank the members of the journal club on infectious diseases at Stockholm University and Daniel Ahlberg for suggestions that lead to substantial improvements of the paper.

\newpage

 \bibliographystyle{plainnat}
\bibliography{main.bib}   


\end{document}